\documentclass[reqno]{amsart}

\usepackage{amssymb, upref, enumerate}

\newtheorem{theorem}{Theorem}[section]

\newtheorem{corollary}[theorem]{Corollary}

\theoremstyle{definition}

\newtheorem{example}[theorem]{Example}

\theoremstyle{remark}

\numberwithin{equation}{section}

\newcommand{\dd}{\mathcal{D}}
\newcommand{\ee}{\mathbb{E}}

\newcommand{\pp}{\mathbb{P}}

\newcommand{\ra}{\rightarrow}
\newcommand{\rr}{\mathbb{R}}

\newcommand{\tr}{\operatorname{Tr}}

\newcommand{\zz}{\mathbb{Z}}
 
\newcommand{\fpar}[2]{\frac{\partial #1}{\partial #2}}

\newcommand{\mpar}[3]{\frac{\partial^2 #1}{\partial #2 \partial #3}}

\newcommand{\ep}{\epsilon}

\newcommand{\mg}{\mathcal{G}}

\newcommand{\mm}{\mathcal{W}}

\newcommand{\mmm}{\widetilde {\mathcal W}}
\newcommand{\wf}{\widetilde{f}}
\newcommand{\WF}{{\widetilde F}}
\newcommand{\wh}{{\widetilde h}}

\newcommand{\hx}{\hat{x}}
\newcommand{\HX}{\hat{X}}
\newcommand{\tilh}{{\widetilde{h}}}
\newcommand\ER{Erd\H{o}s--R\'enyi\ }
\newcommand{\tilf}{{\widetilde{F}}}

\newcommand{\N}{\mathbb{N}}

\begin{document}

\title[Large deviations for random graphs]{An introduction to large deviations for random graphs}


%    author one information
\author{Sourav Chatterjee}
\address{Department of Statistics, Stanford University}
\curraddr{}
\email{souravc@stanford.edu}
\thanks{The author's research was partially supported by NSF grant DMS-1441513}

%    author two information
%\author{}
%\address{}
%\curraddr{}
%\email{}
%\thanks{}

%    \subjclass is required.
%\keywords{Random graph, large deviations, regularity lemma}

\date{}

\dedicatory{}

%    Abstract is required.
\begin{abstract}
This article gives an overview of the emerging literature on large deviations for random graphs. Written for the general mathematical audience, the article begins with a short introduction to the theory of large deviations. This is followed by a description of some large deviation questions about random graphs, and an outline of the recent progress on this topic. A more elaborate discussion follows, with a brief account of graph limit theory and its application in constructing a large deviation theory for dense random graphs. The role of Szemer\'edi's regularity lemma is explained, together with a sketch of the proof of the main large deviation result and some examples. Applications to exponential random graph models are briefly touched upon. The remainder of the paper is devoted to large deviations for sparse graphs. Since the regularity lemma is not applicable in the sparse regime, new tools are needed. Fortunately, there have been several new breakthroughs that managed to achieve the goal by an indirect method. These are discussed, together with an exposition of the underlying theory. The last section contains a list of open problems.
\end{abstract}

\maketitle

\tableofcontents

\section{Large deviations}
The theory of large deviations aims to study two things: (a) the probabilities of rare events, and (b) the conditional probabilities of various events given that some rare event has occurred (that is, what would the world look like if  some rare event happens?).  Often, the second question is more interesting than the first, but it is usually essential to answer the first question to be able to understand how to approach the second. 

By way of illustration, consider the following simple example. Toss a fair coin $n$ times, where $n$ is a large number. Under normal circumstances, you expect to get approximately $n/2$ heads. Also, you expect to get roughly $n/4$ pairs of consecutive heads. However, suppose that the following rare event occurs: the tosses yield $\ge 2n/3$ heads. General-purpose tools from the theory of large deviations allows us to compute that the probability of this rare event is
\begin{equation}\label{ldest}
e^{- n \log (2^{5/3}/3)(1+o(1))}
\end{equation}
as $n\to \infty$.
Moreover, it can be shown that if  this rare event has occurred, then it is highly likely that there are approximately $4n/9$ pairs of consecutive heads instead of the usual $n/4$.

So, how is the estimate \eqref{ldest} obtained? The argument goes roughly as follows. Let $X_1,\ldots, X_n$ be independent random variables, such that $\pp(X_i=0)=\pp(X_i=1)=1/2$ for each $i$. Then the number of heads in $n$ tosses of a fair coin can be modeled by the sum $S_n := X_1+\cdots +X_n$. For any $\theta \ge 0$,
\begin{align*}
\pp(S_n \ge 2n/3) &= \pp(e^{\theta S_n} \ge e^{2\theta n/3})\\
&\le \frac{\ee(e^{\theta S_n})}{e^{2\theta n/3}} \ \ \ \text{(Markov's inequality)}\\
&= \frac{\ee(\prod_{i=1}^n e^{\theta X_i})}{e^{2\theta n /3}}= \frac{\prod_{i=1}^n \ee(e^{\theta X_i})}{e^{2\theta n/3}}\ \ \ \text{(independence)}\\
&= e^{-2\theta n/3} \biggl(\frac{1+e^{\theta}}{2}\biggr)^n\, .
\end{align*}
Optmizing over $\theta$ gives the desired upper bound. To prove the lower bound, take some $\ep>0$ and define a random variable $Z$ as:
\[
Z := 
\begin{cases}
1 &\text{ if } 2n/3 \le S_n \le (1+\ep)2n/3,\\
0 &\text{ otherwise.}
\end{cases}
\]
Then for any $\theta \ge 0$, 
\begin{align*}
\pp(S_n \ge 2n/3) &\ge \ee(Z) \\
&\ge e^{-\theta(1+\ep)2n/3} \ee(e^{\theta S_n} Z) \ \ \ \text{(since $S_n \le (1+\ep)2n/3$ when $Z\ne 0$)}\\
&= e^{-\theta (1+\ep)2n/3} \biggl(\frac{1+e^{\theta}}{2}\biggr)^n\,\frac{\ee(e^{\theta S_n} Z)}{\ee(e^{\theta S_n})}\,.
\end{align*}
The proof is completed by showing that if $\theta$ is chosen to be the same number that optimized the upper bound and $\ep$ is sent to zero sufficiently slowly as $n \to \infty$, then
\[
\frac{\ee(e^{\theta S_n} Z)}{\ee(e^{\theta S_n})} = e^{o(n)}\,.
\]
Establishing the above claim is the most nontrivial part of the whole argument, but is by now standard. This is sometimes called the `change of measure trick'.

The above example has a built-in linearity, which allowed us to explicitly compute $\ee(e^{\theta S_n})$. Generalizing this idea, classical large deviations theory possesses a collection of powerful tools to deal with linear functionals of independent random variables, random vectors, random functions, random probability measures and other abstract random objects. The classic text of Dembo and Zeitouni \cite{dembozeitouni} contains an in-depth  introduction to this broad area.

\section{The problem with nonlinearity}
In spite of the remarkable progress with linear functionals, there are no general tools for large deviations of nonlinear functionals. Nonlinearity arises naturally in many contexts. For instance, the analysis of real-world networks has been one of the most popular scientific endeavors in the last two decades, and  rare events on networks are often nonlinear in nature. This is demonstrated by the following simple example. 

Construct a random graph on $n$ vertices by putting an undirected edge between any two with probability $p$, independently of each other. This is known as the Erd\H{o}s--R\'enyi $G(n,p)$ model, originally defined in \cite{erdosrenyi60}. The model is too simplistic to be a model for any real-world network, but has many nice mathematical properties and has led to the developments of many new techniques in combinatorics and probability theory over the years. One can ask the following large deviation questions about this model: 
\begin{itemize}
\item[(a)] What is the probability that the number of triangles in a $G(n,p)$ random graph is at least $1+\delta$ times the expected value of the number of triangles, where $\delta$ is some given number? 
\item[(b)] What is the most likely structure of the graph, if we know that the above  rare event has occurred?
\end{itemize}
This is an example of a  nonlinear problem, because the number of triangles in $G(n,p)$ is a degree three polynomial of independent random variables. To see this, let $\{1,\ldots,n\}$ be the set of vertices, and let $X_{ij}$ be the random variable that is $1$ if the edge $\{i,j\}$ is present in the graph and $0$ if not. Then $(X_{ij})_{1\le i<j\le n}$ are independent random variables, and the number of triangles is nothing but
\[
\frac{1}{6}\sum_{i,j,k=1}^n X_{ij}X_{jk}X_{ki}\,,
\]
which is a polynomial of degree three. 
Until even a few years ago, large deviations theory did not have the tools to answer such basic questions about nonlinear functions of independent random variables, although a number of powerful concentration inequalities were available for computing upper and lower bounds on tail probabilities~\cite{kimvu00, latala97, talagrand95, vu02}.

\section{Recent developments}\label{prelim}
The large deviation theory for the Erd\H{o}s--R\'enyi random graph was developed fairly recently in~\cite{cv}, taking to completion a program initiated in the unpublished manuscript~\cite{bolthausenetal09}. The theory brought together ideas from classical large deviations theory and tools from combinatorics and graph theory, such as Szemer\'edi's regularity lemma  and the theory of graph limits. The calculations dictated by the theory led to surprising conclusions, even in the simplest of applications such as the following.  Let $T_{n,p}$ be the number of triangles in the Erd\H{o}s--R\'enyi graph $G(n,p)$. What is the most likely structure of the graph if the rare event 
\[
E:=\{T_{n,p} \ge (1+\delta)\ee(T_{n,p})\}
\]
happens, where $\delta$ is a given positive constant? For instance, are all the extra triangles likely to be arising from a small subset of vertices with high connectivity amongst themselves?  Or do they occur because the graph has an excess number of edges spread uniformly? 

Surprisingly, the large deviation theory of \cite{cv} implies that both scenarios can happen. If $p$ is smaller than a threshold, then there exist $0<\delta_1<\delta_2$ such that if $0< \delta \le \delta_1$ or $\delta \ge \delta_2$, then conditional on the event $E$, the graph behaves like $G(n,r)$ for some $r> p$; and  if   $\delta_1< \delta < \delta_2$, then the conditional structure is {\it not} like an Erd\H{o}s--R\'enyi graph. %Explicit formulas for $\delta_1$ and $\delta_2$ were derived by Lubetzky and Zhao \cite{lubetzky12}. 

In other words, if the number triangles exceeds the expected value by a little bit or by a lot, then the most likely scenario is that there is an excess number of edges spread uniformly; and if the surplus amount belongs to a middle range, then the structure of the graph is likely to be inhomogeneous. There is probably no way that the above result could have been guessed from intuition; it was derived purely from a set of mathematical formulas. 

The general theory of \cite{cv} and its main results are described in Section \ref{ldsec}, after a brief introduction to graph limit theory in Section \ref{graphsec}.

The large deviation theory for the Erd\H{o}s--R\'enyi model has been extended to other, more realistic models of random graphs. For example, it was applied to exponential random graph models in \cite{cd3} and a number of subsequent papers. These models are widely used in the analysis of real social networks. A brief discussion of applications to exponential random graphs is given in Section \ref{expsec}. %The theory will be discussed in greater detail in Sections~\ref{graphsec} and~\ref{ldsec}, and its applications to exponential random graph models will be discussed in Section~\ref{ergmsec}. 

In spite of its successes, the theory developed in \cite{cv} has one serious limitation: it applies only to dense graphs. A graph is called dense if the average vertex degree is comparable to the total number of vertices (recall that  the number of neighbors of a vertex is called its degree). For example, in the Erd\H{o}s--R\'enyi model with $n = 10000$ and $p= .3$, the average degree is roughly $3000$.  This is not true for real networks, which are usually sparse.  Unfortunately, the graph theoretic tools used for the analysis of large deviations for random graphs are useful only in the dense setting. In spite of considerable progress in developing a theory of sparse graph limits \cite{bollobasriordan09, bccz1, bccz2}, there is still no result that fully captures the power of Szemer\'edi's lemma in the sparse setting. In the absence of such tools, a nascent theory of `nonlinear large deviations', developed in~\cite{chatterjeedembo14}, has been helpful in solving some questions about large deviations for sparse random graphs~\cite{bglz, lz2}. This theory is discussed in Sections \ref{sparsesec} through \ref{meanfieldsec}.

\section{Graph limit theory}\label{graphsec}
A beautiful unifying theory of graph limits has been developed by Laszlo Lov\'asz and coauthors in recent years~\cite{borgsetal06, borgsetal08, borgsetal12, lovaszbook, lovaszszegedy06}. For connections with the theory of exchangeable arrays in probability theory, see \cite{aldous81, diaconisjanson08, hoover82, kallenberg05}. This section contains a brief review of some of the basic definitions and results from this theory. 

Let $\{G_n\}_{n\ge 1}$ be a sequence of simple graphs whose number of nodes tends to infinity.  For every fixed simple graph $H$, let $\hom(H, G)$ denote the number of homomorphisms of $H$ into $G$ (that is, edge-preserving maps from $V(H)$ into $V(G)$, where $V(H)$ and $V(G)$ are the vertex sets). As an example, note that if $H$ is a triangle, then $\hom(H,G)$ is the number of triangles in $G$ multiplied by~six. 

The number of homomorphisms is normalized to get the homomorphism density 
\begin{equation*}%\label{homdens}
t(H,G) := \frac{\hom(H, G)}{|V(G)|^{|V(H)|}}. 
\end{equation*}
This gives the probability that a random mapping $V(H) \ra V(G)$ is a homomorphism. 

Suppose that $t(H, G_n)$ tends to a limit $t(H)$ for every $H$. Then Lov\'asz and Szegedy \cite{lovaszszegedy06} proved that there is a natural `limit object' in the form of a function $f\in \mm$,  where $\mm$ is the space of all measurable functions from $[0,1]^2$ into $[0,1]$ that satisfy $f(x,y)=f(y,x)$ for all $x,y$. Conversely, every such function arises as the limit of an appropriate graph sequence. This limit object determines all the limits of subgraph densities, as follows. If $H$ is a simple graph on $\{1,2,\ldots,k\}$ with edge set $E(H)$ and $f\in \mm$, let
\begin{equation*}\label{tfdef}
t(H,f) := \int_{[0,1]^k}\prod_{\{i,j\}\in E(H)} f(x_i, x_j) \;dx_1\cdots dx_k. 
\end{equation*}
A sequence of graphs $\{G_n\}_{n\ge 1}$ is said to converge to $f$ if for every finite simple graph $H$,
\begin{equation*}\label{gconv}
\lim_{n\ra \infty} t(H, G_n) = t(H,f).
\end{equation*}
\begin{example}[Limit of Erd\H{o}s--R\'enyi graphs]
Fix $p\in (0,1)$ and let $G_{n,p}$ be a random graph from the Erd\H{o}s-R\'enyi  $G(n,p)$ model. For any fixed graph $H$, it is not difficult to show that with probability $1$,
\[
t(H,G_{n,p}) \ra p^{|E(H)|} \ \text{ as $n\ra\infty$.} 
\]
On the other hand, if $f$ is the function that is identically equal to $p$, then $t(H,f)= p^{|E(H)|}$. 
Thus, the sequence of random graphs $\{G_{n,p}\}_{n\ge 1}$ converges with probability~$1$ to the non-random limit function $f(x,y)\equiv p$ as $n\ra\infty$.
\end{example}
The elements of $\mm$ are sometimes called `graphons'. A finite simple graph $G$ on $\{1,\ldots,n\}$ can be represented as a graphon $f^G$ in a natural way:
\begin{equation*}\label{wg}
f^G(x,y) = 
\begin{cases}
1 &\text{ if $\{[nx], [ny]\}$ is an edge in $G$,}\\
0 &\text{ otherwise.} 
\end{cases}
\end{equation*}
Note that this allows {\it all} simple graphs, irrespective of the number of vertices, to be represented as elements of the single abstract space $\mm$. The starting point of graph limit theory is to define a suitable topology on this space. The first step in defining this topology is to recall the cut distance of Frieze and Kannan \cite{friezekannan99}:
\begin{equation*}\label{defcut}
d_\square (f,g) := \sup_{S,T} \biggl|\int_{S\times T} (f(x,y)-g(x,y)) \,dx \,dy\biggr|\,,
\end{equation*}
where the supremum is taken over all measurable subsets $S$ and $T$ of $[0,1]$. 
The next step is to introduce an equivalence relation on $\mm$, by declaring that $f\sim g$ if  $f(x,y)=g_\sigma(x,y) := g(\sigma x, \sigma y)$ for some measure preserving bijection $\sigma$ of $[0,1]$. Denote by ${\widetilde g}$ the closure  in $(\mm, d_\Box)$ of  the  orbit $\{g_\sigma\}$. Let $\mmm$ be the quotient space and let $\tau$ denote the quotient map $g\mapsto{\widetilde g}$.  Since  $d_\Box$ is invariant under $\sigma$ one can define on $\mmm$ an induced metric $\delta_\Box$ by
$$
\delta_\Box({\widetilde f},{\widetilde g}):=\inf_\sigma d_\Box(f, g_\sigma)=\inf_\sigma d_\Box(f_\sigma, g)=\inf_{\sigma_1,\sigma_2}d_\Box(f_{\sigma_1}, g_{\sigma_2})\,,
$$
making $(\mmm, \delta_\Box)$ a metric space. This is the abstract space of graph limits. To any  finite graph $G$, associate the natural graphon $f^G$ and its orbit  ${\widetilde G}=\tau f^G= {\widetilde f^G}\in\mmm$. One of the key results of the theory is the following:
\begin{theorem}[\cite{borgsetal08}]\label{borgsthm}
A sequence of graphs $\{G_n\}_{n\ge 1}$ converges to a limit $f\in \mm$ if and only if $\delta_\Box({\widetilde G}_n, {\widetilde f}) \ra 0$ as $n \ra \infty$.
\end{theorem} 
Another important result is:
\begin{theorem}[\cite{lovaszszegedy06}]\label{llthm}
The space $\mmm$ is compact under the metric $\delta_\Box$. 
\end{theorem}
The main ingredient in the proofs of the above results is the famous regularity lemma of Szemer\'edi \cite{szemeredi78}. This will be discussed in the next section.

\section{Large deviations for dense random graphs}\label{ldsec}
We will continue to use the notations and terminologies introduced in the previous section. Fix $p\in (0,1)$. 
For $u\in [0,1]$, let 
\begin{equation}\label{ipdef}
I_p(u) := u\log \frac{u}{p}+(1-u)\log \frac{1-u}{1-p}\,,
\end{equation}
with the convention that $0\log 0 = 0$. 
For $h\in \mm$, let 
\begin{equation}\label{iphdef}
I_p(h) := \int_{[0,1]^2} I_p(h(x,y)) \,dx\, dy\,.
\end{equation}
Finally, for $\widetilde{h}\in \mmm$, let $I_p(\widetilde{h}) := I_p(h)$ 
where $h$ is any element of $\widetilde{h}$. A lemma in \cite{cv} shows that the right side does not depend on the choice of $h$ in $\widetilde{h}$, which ensures that the above definition makes sense. 

The  Erd\H{o}s--R\'enyi $G(n,p)$ model induces  a  probability measure $\widetilde{\pp}_{n,p}$  on the space $\mmm$ through the map $G\ra \widetilde{G}$. The main result of \cite{cv} is the following large deviation principle for $\widetilde{\pp}_{n,p}$.
\begin{theorem}[\cite{cv}]\label{cvthm}
For any closed set $\widetilde {F} \subseteq\mmm$,
\begin{align*}
\limsup_{n\ra\infty} \frac{2}{n^2}\log {\widetilde\pp}_{n,p}(\widetilde {F}) &\le -\inf_{{\widetilde h}\in \widetilde{F}} I_p({\widetilde h})\,,
\end{align*}
and for any open set $\widetilde{U}\subseteq \mmm$,
\begin{align*}
 \liminf_{n\ra\infty} \frac{2}{n^2}\log {\widetilde\pp}_{n,p}( \widetilde{U}) &\ge -\inf_{{\widetilde h}\in \widetilde{U}} I_p({\widetilde h})\,. 
\end{align*}
\end{theorem}
Although the above theorem looks like a result specifically about the $G(n,p)$ model, this is somewhat misleading. Theorem \ref{cvthm} actually allows us to approximately count the number of simple graphs on $n$ vertices that have any given property, as long as the property is nicely behaved with respect to the cut metric. This can be made precise as follows. For any Borel set $\widetilde{A}\subseteq \mmm$, let 
\[
\widetilde{A}_n := \{\wh\in \widetilde{A}: \wh = \widetilde{G} \text{ for some $G$ on $n$ vertices}\}\,. 
\]
Let 
\[
I(u) := u\log u + (1-u)\log (1-u)\,.
\]
For any $\widetilde{h}\in \mmm$, let 
\begin{equation}\label{idef}
I(\widetilde{h}) := \int_{[0,1]^2} I(h(x,y)) \,dx\, dy\,,
\end{equation}
where $h$ is any element of $\widetilde{h}$. The following corollary can be easily derived from Theorem \ref{cvthm}, by taking $p=1/2$.
\begin{corollary}\label{cvcor}
For any measurable $\widetilde{A}\subseteq \mmm$,
\begin{align*}
-\inf_{\wh\in \textup{cl}(\widetilde{A})} I(\wh)&\ge \limsup_{n\ra\infty}\frac{2\log|\widetilde{A}_n|}{n^2}\\
&\ge \liminf_{n\ra\infty}\frac{2\log|\widetilde{A}_n|}{n^2} \ge -\inf_{\wh\in \textup{int}(\widetilde{A})}I(\wh)\,,
\end{align*}
where $\textup{cl}(\widetilde{A})$ is the closure of $\widetilde{A}$ and $\textup{int}(\widetilde{A})$ is the interior of $\widetilde{A}$.
\end{corollary}
Under very special circumstances, the variational problems of Theorem \ref{cvthm} and Corollary \ref{cvcor} are known to have an explicit solutions. For example, it follows from Corollary \ref{cvcor} and some results in \cite{cv} that the number of graphs on $n$ vertices with at least $tn^3$ triangles is 
\[
e^{\frac{1}{2}n^2 f(t)(1+o(1))}
\]
as $n\to \infty$, where 
\[
f(t) =
\begin{cases}
\log 2 &\text{ if } 0\le t< \frac{1}{48},\\
-I((6t)^{1/3}) &\text{ if } \frac{1}{48} \le t \le \frac{1}{6},\\
-\infty &\text{ if } t > \frac{1}{6}. 
\end{cases}
\]
On the other hand, for the number of graphs with {\it at most} $tn^3$ triangles, such an explicit formula can be obtained if $t$ is sufficiently away from zero, and it can also be shown that this formula {\it does not} hold if $t$ is sufficiently close to zero. As of now, there is no explicit formula for small $t$. See Zhao~\cite{zhao15} for the most advanced results about the lower tail problem.

Theorem \ref{cvthm} gives estimates of the probabilities of rare events related to an Erd\H{o}s--R\'enyi graph. However, it does not answer the second type of question mentioned in the introduction, that is, given that some particular rare event has occurred, what does the graph look like? The answer may be roughly described as follows. If $G_{n,p}$ is a random graph from the $G(n,p)$ model,  $\WF$ is a closed subset of $\mmm$, and we condition on the event that $\widetilde{G}_{n,p}\in \WF$, then it is very likely that $\widetilde{G}_{n,p}$ is close to one of the elements of $\WF$ that minimize $I_p$ in $\WF$. The precise result goes as follows. Let $\WF$ and $G_{n,p}$ be as above. Let $\WF^*$ be the set of minimizers of $I_p$ in $\WF$. Suppose that 
\begin{equation}\label{inf0}
\inf_{\wh\in \textup{int}(\WF)} I_p(\wh) = \inf_{\wh \in \WF} I_p(\wh)> 0\,.
\end{equation}
It was proved in \cite{cv} that $I_p$ is lower semicontinuous on $\mmm$, and we know that $\mmm$ is compact by Theorem \ref{llthm}. Therefore $\WF^*$ is nonempty. For any $\wh\in \mmm$, let
\begin{equation}\label{deltagf}
\delta_\Box(\wh,\WF^*) := \inf_{{\widetilde f}\in \WF^*} \delta_\Box(\wh, \wf)\,. 
\end{equation}
\begin{theorem}[\cite{cv}]\label{conditional}
In the above setting, the following inequality holds for each $n\ge 1$ and $\ep >0$:
\[
\pp(\delta_\Box(\widetilde{G}_{n,p}, \WF^*) \ge \ep \mid \widetilde{G}_{n,p}\in \WF) \le e^{-C(\ep, \WF) n^2}
\]
where $C(\ep, \WF)$ is a positive constant depending only on $\ep$ and $\WF$.  
\end{theorem}
In particular, it follows that if $\WF^*$ contains only one element, $\wh^*$, then the conditional distribution of $\widetilde{G}_{n,p}$ given $\widetilde{G}_{n,p}\in \WF$ converges to the point mass at $\wh^*$ as $n\ra\infty$, giving a conditional law of large numbers. 

Let us now see a sketch of the proof of Theorem \ref{cvthm}. The result can be proved by standard techniques for the weak topology on $\mmm$. However, the weak topology is not very interesting. For example, subgraph counts are not continuous with respect to the weak topology. The large deviation principle under the topology of the cut metric (Theorem \ref{cvthm}) does not follow via standard methods.

The main tool for proving Theorem \ref{cvthm} is Szemer\'edi's regularity lemma. One version of Szemer\'edi's lemma goes as follows. Let $G= (V,E)$ be a simple graph of order $n$ (recall that the number of vertices of a graph is called its order). For any $X,Y\subseteq V$, let $e_G(X,Y)$ be the number of $X$-$Y$ edges of $G$ and let
\[
\rho_G(X,Y) := \frac{e_G(X,Y)}{|X||Y|}\,,
\]
where $|X|$ and $|Y|$ are the sizes of $X$ and $Y$.  Call a pair $(A,B)$ of disjoint sets $A, B \subseteq V$ $\ep$-regular if all $X\subseteq A$ and $Y\subseteq B$ with $|X|\ge \ep|A|$ and $|Y|\ge \ep |B|$ satisfy 
\[
|\rho_G(X,Y) - \rho_G(A,B)| \le \ep\,. 
\]
The concept of $\ep$-regularity tries to capture the notion that the edges going between $A$ and $B$ behave like randomly distributed edges with density $\rho_G(A,B)$.

A partition $\{V_0,\ldots,V_K\}$ of $V$ is called an $\ep$-regular partition of $G$ if it satisfies the following conditions: 
\begin{itemize}
\item[(i)] $|V_0|\le \epsilon n$.
\item[(ii)] $|V_1|=|V_2|=\cdots= |V_K|$.
\item[(iii)]  All but at most $\ep K^2$ of the pairs $(V_i, V_j)$, $1\le i<j\le K$, are $\ep$-regular. 
\end{itemize}
\begin{theorem}[Szemer\'edi's regularity lemma \cite{szemeredi78}]\label{szthm}
Given any $\ep >0$ and $m \ge 1$ there exists $M = M(\ep,m)$ such that every graph  of order $\ge M$ admits an $\ep$-regular partition $\{V_0,\ldots, V_K\}$ for some $K\in [m, M]$. 
\end{theorem}
Roughly speaking, Szemer\'edi's lemma says that any large graph $G$ may be partitioned into blocks of equal size (plus one exceptional block of small size) so that the edges going between the blocks behave like randomly distributed edges, and the number of blocks depends only on the desired degree of randomness and not on the size of the graph. The key to proving Theorem \ref{cvthm} is to formulate a precise version of this statement. This goes as follows.

Choose a small number $\ep>0$. Suppose that $G = (V,E)$ is a simple graph of order $n$ with $\ep$-regular partition $\{V_0,\ldots, V_K\}$, as in Theorem \ref{szthm}. Let $G' = (V, E')$ be a random graph with independent edges where a vertex $u\in V_i$ is connected to a vertex $v\in V_j$ with probability $\rho_G(V_i, V_j)$. Let $\widetilde{\mu}$ be the probability measure on $\mmm$ induced by $G'$. The main step in the proof of Theorem \ref{cvthm}, proved using Szemer\'edi's lemma, is that $\widetilde{\mu}(\widetilde{B}) \approx 1$, where $\widetilde{B}$ is a small ball around $\widetilde{G}$ in the $\delta_\Box$ metric and the radius of $\widetilde{B}$  depends only on $\ep$ and not on the size of $G$. In other words, the random graph $G'$ is close to the given graph $G$ in the $\delta_\Box$ metric. This gives a precise meaning to the sentence from the previous paragraph that a large graph can be partitioned into blocks with approximately randomly distributed edges between the blocks. 

To complete the sketch of the proof of Theorem \ref{cvthm}, let $f$ be the probability density of $\widetilde{\pp}_{n,p}$ with respect to $\widetilde{\mu}$. Since the edges in both $G(n,p)$ and $G'$ are independent, this probability density is easy to compute.  Since $\widetilde{B}$ is a ball of small radius, we get
\[
\widetilde{\pp}_{n,p}(\widetilde{B}) \approx f(\widetilde{G}) \widetilde{\mu}(\widetilde{B}) \approx f(\widetilde{G})\,. 
\]
Since the space $\mmm$ is compact (Theorem \ref{llthm}), this allows us to estimate $\widetilde{\pp}(\widetilde{A})$ for any nice  set $\widetilde{A}$ by approximating $\widetilde{A}$ as a finite union of small balls. 
\begin{example}[Triangles in dense Erd\H{o}s--R\'enyi graphs]
For an application of Theorem \ref{cvthm} to a concrete problem, consider the number of triangles $T_{n,p}$ in a $G(n,p)$ random graph. Recall that $\mm$ is the space of symmetric measurable functions from $[0,1]^2$ into $[0,1]$.  For each $f\in \mm$, let
\[
T(f) := \frac{1}{6} \int_{[0,1]^3} f(x,y) f(y,z) f(z,x) \,dx\,dy\,dz\,,
\]
and let $I_p(f)$ be as in \eqref{iphdef}. For each $p\in (0,1)$ and $t\ge 0$, let
\begin{equation}\label{phidef}
\phi_p(t) := \inf\{I_p(f) : f\in \mm,\,T(f) \ge t\}. 
\end{equation}
In \cite{cv}, the following result was proved using Theorem \ref{cvthm}.
\begin{theorem}[\cite{cv}]\label{triangle}
For each $p\in (0,1)$ and each $t\ge 0$, 
\begin{align*}
\lim_{n\ra\infty} \frac{2}{n^2}\log \pp(T_{n,p}\ge tn^3) =  -\phi_p(t). 
\end{align*}
Moreover, the infimum is attained in the variational problem \eqref{phidef}. 
\end{theorem}
The way to deduce Theorem \ref{triangle} from Theorem \ref{cvthm} is to transform the problem into a question about the probability measure $\widetilde{\pp}_{n,p}$. This is quite straightforward. The next step is to apply Theorem \ref{cvthm} to an appropriate pair of sets $(\widetilde{F},\widetilde{U})$ and show that the upper and lower bounds match. 

As a consequence of Theorem \ref{llthm}, one can show that the infimum in \eqref{phidef} is attained. Moreover, the minimizing functions determine the behavior of the graph conditional on the event $\{T_{n,p}\ge tn^3\}$. Theorem \ref{borgsthm} is instrumental in proving such claims. For example, if the minimization problem is solved by constant functions, then the conditional behavior continues to be like that of an Erd\H{o}s--R\'enyi graph, but with a different edge probability. On the other hand, if the minimizers are all non-constant functions, then the conditional behavior is not like that of an Erd\H{o}s--R\'enyi graph. The precise statements and proofs of these claims can be formalized using Theorem \ref{conditional}.

The variational problem \eqref{phidef} was analyzed in \cite{cv} and \cite{lubetzky12}. One of the most surprising findings, already mentioned briefly in Section \ref{prelim}, is that if $p$ is smaller than a threshold, then there exists a non-empty interval $(t_1,t_2)$ such that the minimization problem is solved by constant functions when $t\not\in (t_1,t_2)$ and is solved by non-constant functions when $t\in (t_1,t_2)$. The existence of this interval was proved in \cite{cv} and the values of $t_1$ and $t_2$ were computed in \cite{lubetzky12}. The theorem proved in \cite{lubetzky12} applies to the homomorphism density of any regular graph; specialized to triangles, it says the following.
\begin{theorem}[\cite{lubetzky12}]
Take any $p\in (0,1)$ and $t\in (p^3/6,1)$. Let $r:= (6t)^{1/3}$. If the point $(r^2, I_p(r))$ lies on the convex minorant of the function $J_p(x):= I_p(x^{1/2})$, then there is a unique solution of the variational problem \eqref{phidef}, and it is the constant function $f(x,y)\equiv r$. On the other hand, if $(r^2, I_p(r))$ does not lie on the convex minorant of $J_p$, then any solution of \eqref{phidef} is non-constant.
\end{theorem}
\end{example}

Applications of Theorem \ref{cvthm} to exponential random graphs are briefly discussed in the next section, before moving on to large deviations for sparse random graphs. 

\section{Exponential random graphs}\label{expsec}

Let $\mg_n$ be the set of all simple graphs on $n$ labeled
vertices. A variety of probability models on this set can be presented in exponential form
\begin{equation*}
p_\beta(G)=\exp\biggl(\sum_{i=1}^k\beta_iT_i(G)-\psi(\beta)\biggr)
\label{21}
\end{equation*}
where $\beta=(\beta_1,\dots,\beta_k)$ is a vector of real parameters,
$T_1,T_2,\dots,T_k$ are real-valued functions on $\mg_n$, and
$\psi(\beta)$ is the normalizing constant. Usually, $T_i$ are taken to
be counts of various subgraphs, for example $T_1(G)=$ number of edges in $G$,
$T_2(G)=$ number of triangles in $G$, etc. These are known as exponential random graph models (ERGM). These models are widely used in the study of social networks, but were generally mathematically intractable until recently. A general technique for proving theorems about dense exponential random graphs, based on the large deviation theory for dense Erd\H{o}s--R\'enyi random graphs, was proposed in \cite{cd3}. The mathematical literature on exponential random graphs has grown considerably since the publication of \cite{cd3}. Since this is somewhat disjoint from the large deviations literature, I will not attempt to survey these developments here. For a comprehensive survey of the ERGM literature till the publication of \cite{cd3}, see \cite{cd3}. For a few mathematical results preceding \cite{cd3}, see \cite{bbs, chatterjeedey10}. For a non-exhaustive list of subsequent developments, see \cite{ar13, az15,kenyon14, ky14, lubetzky12, rrs,  radinsadun13, radinsadun15, radinyin13, sr13, yin13}. The discussion in this section will be limited to a basic result from \cite{cd3} and one easy example. We will continue to use the notations introduced in Section \ref{graphsec}.

Let $T$ be  a real-valued continuous function on the space $\mmm$. Define a probability mass
function $p_n$ on $\mg_n$ induced by $T$ as:
\begin{equation*}
p_n(G):=e^{n^2 (T(\widetilde{G})-\psi_n)}\,, 
\end{equation*}
where $\psi_n$ is a constant such that the total
mass of $p_n$ is $1$. Explicitly,
\begin{equation*}%\label{psidef}
\psi_n=\frac1{n^2}\log\sum_{G\in\mg_n}e^{n^2 T(\widetilde{G})} 
\end{equation*}
The coefficient $n^2$ is meant to ensure that $\psi_n$ tends to a
non-trivial limit as $n$ tends to infinity. This setup gives an abstract formulation of exponential random graphs in the language of graph limits, but with the usual limitation that it makes sense only for dense graphs. The following theorem was proved in \cite{cd3}.
\begin{theorem}[\cite{cd3}]\label{soln}
Let $T$  and $\psi_n$ be as above. Let $I$ be the function defined in equation \eqref{idef}. Then
\begin{equation}\label{expvar}
\lim_{n\to\infty} \psi_n = \sup_{\widetilde{h}\in\mmm} \Big(T(\widetilde{h}) - \frac{1}{2}I(\widetilde{h})\Big)\,.
\end{equation}
\end{theorem}
The evaluation of the normalizing constant is an important problem in statistical applications of exponential random graphs, because it is required for computing maximum likelihood estimates. Incidentally, even the existence of the limit in Theorem~\ref{soln} has an important consequence. Suppose that a computer program can evaluate the exact value of the normalizing constant for moderate sized $n$. Then if $n$ is large, one can choose a scaled down model with a smaller number of nodes, and use the exact value of the normalizing constant in the scaled down model as an approximation to the normalizing constant in the larger model. 

Theorem~\ref{soln} gives an asymptotic formula for $\psi_n$. However,
it says nothing about the behavior of a random graph drawn from the
exponential random graph model. Some aspects of this behavior can be
described as follows. Let $\tilf^*$ be the subset of $\mmm$ where
$T(\tilh)-\frac{1}{2}I(\tilh)$ is maximized. By the compactness of $\mmm$, the
continuity of $T$ and the lower semi-continuity of $I$ (proved in \cite{cv}), $\tilf^*$ is a
non-empty compact set. Let $G_n$ be a random graph on $n$ vertices
drawn from the exponential random graph model defined by $T$. The
following theorem shows that for $n$ large, $\widetilde{G}_n$ must lie
close to $\tilf^*$ with high probability. In particular, if $\tilf^*$ is a
singleton set, then the theorem gives a weak law of large numbers for
$\widetilde{G}_n$.
\begin{theorem}[\cite{cd3}]\label{limitbehave}
  Let $\tilf^*$ and $G_n$ be defined as in the above paragraph.   Then for
  any $\eta > 0$ there exist $C, \gamma >0$ such that for any $n$,
\begin{equation*}
\pp(\delta_\Box(\widetilde{G}_n, \tilf^*) > \eta) \le Ce^{-n^2  \gamma}\,. 
\end{equation*}
\end{theorem}
There is no general method for efficiently solving the variational problem \eqref{expvar}, either analytically or computationally. In some special cases, however, the problem yields an explicit solution. One such example, worked out in \cite{cd3}, is the following.
\begin{example}
Let $H_1,\dots, H_k$ be finite simple graphs, where $H_1$ is the
complete graph on two vertices (that is, just a single edge), and each
$H_i$ contains at least one edge. Let $\beta_1,\dots,\beta_k$ be $k$
real numbers. For any $h\in \mm$, let 
\begin{equation}\label{tdef1}
T(h) := \sum_{i=1}^k \beta_i t(H_i, h)
\end{equation}
where $t(H_i, h)$ is the homomorphism density of $H_i$ in $h$. The functional $T$ extends naturally to $\mmm$, and is continuous on $\mmm$ by Theorem \ref{borgsthm}. For any finite simple graph
$G$ that has at least as many nodes as the largest of the $H_i$'s,
\begin{equation*}
T(\widetilde{G}) = \sum_{i=1}^k \beta_i t(H_i, G)\,,
\end{equation*}
where $t(H_i, G)$ is the homomorphism density of $H_i$ in $G$. 
For example, if $k=2$, $H_2$ is a triangle and $G$ has at least
three nodes, then
\begin{equation*}
T(\widetilde{G}) = 2\beta_1 \frac{\text{number of edges in } G}{n^2} + 6\beta_2 \frac{\text{number of triangles in } G}{n^3}\,. 
\end{equation*}
The following theorem says that when $T$ is of the form \eqref{tdef1} and  $\beta_2,\dots,\beta_k$ are
nonnegative, the variational problem of Theorem \ref{soln} can be reduced to a simple maximization problem in one real variable. The theorem moreover says that each solution of the variational problem is a constant function, and there are only a finite number of solutions. By Theorem \ref{limitbehave}, this implies that when $\beta_2,\ldots, \beta_k$ are nonnegative, exponential random graphs from this class of models behave like random graphs drawn from a finite mixture of \ER models.
\begin{theorem}[\cite{cd3}]\label{specialnorm}
  Let $H_1,\dots, H_k$ and $T$ be as above. Suppose
  that the parameters $\beta_2,\dots,\beta_k$ are nonnegative. Let $\psi_n$ be the normalizing constant of the exponential random graph model induced by $T$ on the set of simple graphs on $n$ vertices. Then
\begin{equation}\label{scalar}
\lim_{n\to\infty} \psi_n = \sup_{0\le u\le 1} \biggl(\sum_{i=1}^k \beta_i u^{e(H_i)} - \frac{1}{2}I(u)\biggr) 
\end{equation}
where $I$ is the function defined in \eqref{idef} and $e(H_i)$ is the
number of edges in $H_i$.  Moreover, there are only a finite number of solutions of the variational
problem of Theorem~\ref{soln} for this $T$, and each solution is a constant function,
where the constant solves the scalar maximization problem~\eqref{scalar}.
\end{theorem}
\end{example}

\section{Large deviations for sparse random graphs}\label{sparsesec}
Recall the Erd\H{o}s--R\'enyi model $G(n,p)$. Let $T_{n,p}$ be the number of triangles in a $G(n,p)$ random graph, as before. The behavior of the upper tail probabilities of $T_{n,p}$, especially when $p\to 0$ as $n\to \infty$, has been the subject of intense investigation for many years. After a series of successively improving suboptimal results by various authors, a big advance was made by Kim and Vu~\cite{kimvu04} and simultaneously by Janson et al.~\cite{JOR04} who showed that if $p\ge n^{-1}\log n$,  then 
\[
e^{-c_1(\delta) n^2 p^2\log(1/p)}\le \pp(T_{n,p}\ge (1+\delta)\ee(T_{n,p}))\le e^{-c_2(\delta)n^2p^2}\, ,
\]
where $c_1(\delta)$ and $c_2(\delta)$ are constants that depend only on $\delta$.
It took several more years to remove the logarithmic discrepancy between the exponents on the two sides. It was finally established in~\cite{chatriangle} and independently in~\cite{dk, dk2} that when $p\ge n^{-1}\log n$,
\begin{align*}
e^{-c_1(\delta) n^2 p^2\log(1/p)}&\le \pp(T_{n,p}\ge (1+\delta)\ee(T_{n,p}))\le e^{-c_2(\delta)n^2p^2\log(1/p)}\, . 
\end{align*}
This still left open the question of determining the dependence of the exponent on~$\delta$. We have seen the solution of this question in the previous section when $p$ is fixed. But this theory does not carry over to the case where $p\to 0$ as $n\to \infty$, especially if $p$ decays like a negative power of $n$. Partly motivated by this question, a general theory of nonlinear large deviations was proposed  in \cite{chatterjeedembo14}. Using this theory, Lubetzky and Zhao \cite{lz2} proved the following theorem, which fully solved the upper tail large deviation problem for $T_{n,p}$ when $p$ goes to zero slower than~$n^{-1/42}$.
\begin{theorem}[\cite{lz2}]\label{lzthm}
If $T_{n,p}$ is the number of triangles in $G(n,p)$, then as $n\ra\infty$ and $p\ra 0$ slower than $n^{-1/42}$, 
\begin{align*}
&\pp(T_{n,p}\ge (1+\delta)\ee(T_{n,p})) = \exp\biggl(-(1+o(1)) \min\biggl\{\frac{\delta^{2/3}}{2}\, , \, \frac{\delta}{3}\biggr\}n^2 p^2 \log \frac{1}{p}\biggr)\, .
\end{align*}
\end{theorem}
It is conjectured that this result holds when $p\ra 0$ slower than $n^{-1/2}$ (see \cite{lz2}). The above theorem has been generalized by Bhattacharya et al.~\cite{bglz}, who got the following beautiful result. Take any finite simple graph $H$ with maximum degree~$\Delta$. Let $H^*$ be the induced subgraph of $H$ on all vertices whose degree in $H$ is $\Delta$. Recall that an independent set in a graph is a set of vertices such that no two are connected by an edge. Also, recall that a graph is called regular if all its vertices have the same degree, and irregular otherwise. Define a polynomial 
\[
P_{H^*}(x) := \sum_k i_{H^*}(k) x^k\,,
\] 
where $i_{H^*}(k)$ is the number of $k$-element independent sets in $H^*$. The main result of \cite{bglz} is the following. %Theorem \ref{lzthm} is a special case of this result. 
\begin{theorem}[\cite{bglz}]\label{bthm}
Let $H$ be a connected finite simple graph on $k$ vertices with maximum degree $\Delta \ge 2$. Then for any $\delta >0$, there is a unique positive number $\theta = \theta(H, \delta)$ that solves $P_{H^*}(\theta)=1+\delta$, where $P_{H^*}$ is the polynomial defined above. Let $H_{n,p}$ be the number of homomorphisms (defined in Section \ref{graphsec}) of $H$ into an Erd\H{o}s--R\'enyi $G(n,p)$ random graph. Then there is a constant $\alpha_H>0$ depending only on $H$, such that if $n\to \infty$ and $p\to 0$ slower than $n^{-\alpha_H}$, then for any $\delta >0$,
\[
\pp(H_{n,p} \ge (1+\delta) \ee(H_{n,p})) = \exp\biggl(-(1+o(1))c(\delta) n^2 p^\Delta \log \frac{1}{p}\biggr)\,,
\]
where 
\[
c(\delta) = 
\begin{cases}
\min\{\theta,\frac{1}{2}\delta^{2/k}\} & \text{ if $H$ is regular,} \\
\theta &\text{ if $H$ is irregular.}
\end{cases}
\]
%Here $|V(H)|$ and $|E(H)|$ stand for the number of vertices and the number of edges in $H$. 
\end{theorem}
The formula given in Theorem \ref{bthm} is more than just a formula. It gives a hint at the conditional structure of the graph, and at the nature of phase transitions as $\delta$ varies. Unlike the dense case, it is hard to give a precise meaning to claims about the conditional structure in the sparse setting due to the lack of an adequate sparse graph limit theory. For a detailed discussion, see \cite{bglz}.

The paper \cite{bglz} also gives a number of examples where the coefficient $c(\delta)$ in Theorem \ref{bthm} can be explicitly computed. For instance, if $H = C_4$, the cycle of length four, then
\[
c(\delta) = 
\begin{cases}
\frac{1}{2}\sqrt{\delta} &\text{ if } \delta<16,\\
-1+\sqrt{1+\frac{1}{2}\delta} &\text{ if } \delta \ge 16.
\end{cases}
\]
Theorem \ref{lzthm} is also a special case of Theorem \ref{bthm}.

\section{The low complexity gradient condition}\label{lowsec}
The purpose of this section is to begin to describe the theory developed in \cite{chatterjeedembo14} that leads to the results of previous section. The initial part of the discussion will be kept intentionally imprecise so as to convey the ideas smoothly without getting bogged down in technical details and notational complexities. The precise statement of the main result is given in  in the next section.

Take a smooth $f:[0,1]^n \ra \rr$. Let $Y = (Y_1,\ldots, Y_n)$ be a vector of independent $0$-$1$ random variables with $\pp(Y_i=1)=p$. The goal is to find an approximation for the upper tail probability $\pp(f(Y)\ge tn)$ when $t$ is bigger than $n^{-1}\ee(f(Y))$. Recall the function $I_p:[0,1]\to \rr$ defined in \eqref{ipdef}. For $x= (x_1,\ldots, x_n)$, define 
\begin{equation*}\label{ipdef2}
I_p(x) := \sum_{i=1}^nI_p(x_i)\, .
\end{equation*}
For each $t\in \rr$, let
\begin{equation}\label{phipdef}
\phi_p(t) := \frac{1}{n}\inf\{I_p(x)\, : \, x\in [0,1]^n,\, f(x)\ge tn\}\, . 
\end{equation}
In many problems, it turns out that
\begin{equation}\label{star}
\pp(f(Y) \ge t n) \approx e^{-\phi_p(t)n}\, .
\end{equation}
%Note that the left side depends only the values of $f$ on $\{0,1\}^n$, whereas the right side depends on the values of $f$ on the full domain $[0,1]^n$.
The approximation \eqref{star} holds in great generality for linear functions. Theorem~\ref{triangle} gives a nonlinear example. Theorems~\ref{lzthm} and~\ref{bthm} are consequences of this approximation where $\phi_p(t)$ is explicitly computable in a limit.  The main result of~\cite{chatterjeedembo14}, building on ideas developed in~\cite{chathesis,cha07,chatterjeedey10}, gives a sufficient condition under which~\eqref{star} is valid for a general nonlinear $f$. The condition may be roughly stated as follows. Let $\nabla f$ be the gradient of $f$, so that $\nabla f$ is a map from $[0,1]^n$ into $\rr^n$. Then we need that the image of $[0,1]^n$ under the map $\nabla f$ has {\it low complexity}, in the sense that it can be covered by a relatively small number of balls of an appropriate radius. Here `small' means $e^{o(n)}$. This is called the `low complexity gradient condition' in \cite{chatterjeedembo14}. A different way of putting this is to say that the value of $\nabla f(x)$ may be approximately encoded by $o(n)$ bits of information.

To understand this, let $f(x) = \frac{1}{2}\|x\|^2$, where $\|x\|$ is the Euclidean norm of $x$. Then $\nabla f(x)=x$, and so $\nabla f([0,1]^n)=[0,1]^n$. Since the typical distance between points in $[0,1]^n$ is of order $\sqrt{n}$, it is appropriate to consider balls whose radii are of order $\sqrt{n}$ when measuring the complexity on $[0,1]^n$. An $\ep\sqrt{n}$-net of $[0,1]^n$ is a collection of points in $[0,1]^n$ such that any point of $[0,1]^n$ is within distance $\ep\sqrt{n}$ from one of these points. It is not difficult to prove that for any fixed $\ep$, an $\ep\sqrt{n}$-net of $[0,1]^n$ must have size at least $e^{C(\ep)n}$, where $C(\ep)$ is a positive constant that depends on $\ep$. Therefore the gradient of $f$ does not have low complexity. 

On the other hand, suppose that $f:[0,1]^n \ra\rr$ is a linear map. Then $\nabla f(x)$ is the same for every $x$, and hence the image of $[0,1]^n$ under $\nabla f$ consists of a single point. This is a set of the lowest possible complexity, and therefore the low complexity gradient condition is satisfied by linear functions. We will see  less trivial examples of functions satisfying the low complexity gradient condition below, after formally defining a measure of the complexity of $\nabla f$. This measure of complexity is chiefly for expositional purposes; the main theorem of \cite{chatterjeedembo14}, presented later in this section, does not make any direct use of this complexity measure.

Implicitly, we have a sequence of functions rather than a single function; that is, the function $f$ in \eqref{star} depends on $n$ and varies as $n$ varies. We will assume that $f$ is scaled with $n$ in such a way that sizes of the partial derivatives of $f$ do not go to zero or infinity as $n\ra\infty$. Consequently, the typical diameter of $\nabla f([0,1]^n)$ is of order $\sqrt{n}$. Keeping this in mind, the complexity of $\nabla f([0,1]^n)$ may be defined as follows. For each $\ep$, let $\dd(\ep)$ be an $\ep\sqrt{n}$-net for the image of $[0,1]^n$ under $\nabla f$. As explained before, this is a finite set of points such that for any $x\in [0,1]^n$, there exists a point $z\in \dd(\ep)$ such that the Euclidean distance between $\nabla f(x)$ and $z$ is less than $\ep\sqrt{n}$. Let $|\dd(\ep)|$ denote the size of $\dd(\ep)$. Suppose that $\dd(\ep)$ is optimized to have the minimum possible size. Define
\[
C(f) := \inf_{\ep>0}\biggl(\ep + \frac{\log |\dd(\ep)|}{n}\biggr)\,.
\]
Note that $C(f)$ is small if and only if there exists a small $\ep$ such that $\log |\dd(\ep)| \ll n$. The number $C(f)$ can therefore be used as a measure of the complexity of the gradient of $f$ when $f$ scales with $n$ is such a way that the sizes of its partial derivatives do not tend to zero or infinity as $n\ra\infty$. We have already observed that $C(f)=0$ when $f$ is linear. Let us now see some more examples.
\begin{example}
For $x=(x_1,\ldots, x_n)\in [0,1]^n$, let 
\[
f(x):=\frac{1}{n}\sum_{1\le i<j\le n} x_i x_j\, .
\]
Then for each $i$,
\[
\fpar{f}{x_i}= \frac{1}{n}\sum_{j\ne i} x_j = -\frac{x_i}{n} + \frac{1}{n}\sum_{j=1}^n x_j\, .
\]
Note that the sizes of the partial derivatives of $f$ are not blowing up or tending to zero as $n\ra \infty$. Let $D$ be the set of all vectors in $\rr^n$ that are of the form $(k/n, k/n, \ldots, k/n)$ for some  $0\le k\le n-1$. The above formula shows that for any $x\in [0,1]^n$, there exists $z\in D$ such that the Euclidean distance between $\nabla f(x)$ and $z$ is bounded by $2n^{-1/2}$. To see this, just take $z=(k/n,\ldots, k/n)$ where $k$ is the integer part of $\sum_{j} x_j$. Thus, the set $D$ can serve as $\dd(\ep)$ for this $f$, for $\ep = 2n^{-1}$. Since $|\dd(\ep)|=n$, we get 
\[
C(f)\le \frac{2+\log n}{n}\,,
\]
which shows that $f$ satisfies the low complexity gradient condition. This is easier to understand via the language of encoding $\nabla f(x)$: the value of $\nabla f(x)$ may be approximately encoded by the single quantity $n^{-1}\sum_j x_j$, and therefore needs $O(\log n)$ bits.
\end{example}
\begin{example}%[The Hamiltonian of the 1D Ising model]
Take any $n\ge 2$ and for $x= (x_1,\ldots,x_n)\in [0,1]^n$, let 
\[
f(x) := \sum_{i=1}^{n-1} x_i x_{i+1}\, .
\]
Then
\[
\fpar{f}{x_i} = 
\begin{cases}
x_2 &\text{ if } i=1,\\
x_{i-1} + x_{i+1} &\text{ if } 2\le i\le n-1,\\
x_{n-1} &\text{ if } i=n.
\end{cases}
\]
This $f$, again, is scaled in such a way that the sizes of its partial derivatives remain stable as $n\ra\infty$. Now, given any $z\in \{0,1,2\}^n$, it is easy to prove using the above formula that there exists $x\in \{0,1\}^n$ such that $\nabla f(x)$ and $z$ agree on at least $n/3$ coordinates. This, in turn, can be used to show  that $C(f)$ does not tend to zero as $n\ra\infty$.  Therefore, this $f$ does not have a gradient of low complexity. Again, this may be easier to understand by observing that to encode $\nabla f(x)$, we need to essentially have information about all the coordinates of $x$, which requires $n$ bits.%One can check that the approximation \eqref{star} is not valid for this~$f$.
\end{example}

\begin{example}[Triangles in $G(n,p)$]\label{tex}
Take any $n$ such that $n = m(m-1)/2$ for some positive integer $m\ge 2$. Denote elements of $\rr^n$ as $x = (x_{ij})_{1\le i<j\le m}$. Define  $f:\rr^n \ra \rr$ as 
\[
f(x)=\frac{1}{m}\sum_{i,j,k=1}^mx_{ij}x_{jk}x_{ki}\, ,
\]
using the conventions that $x_{ii}=0$ for each $i$ and $x_{ij}=x_{ji}$ when $i>j$. Note that when $x_{ij}$ is the $(i,j)^{\textup{th}}$ entry of the adjacency matrix of a graph on $m$ vertices, $f(x)$ is equal to the number of triangles in the graph multiplied by $6/m$. Let us now see why $C(f)$ is small. The complete proof  is somewhat lengthy, so we will only see a sketch of the proof. For details, see \cite{chatterjeedembo14}. 

First, a simple computation gives
\begin{equation}\label{tder}
\fpar{f}{x_{ij}} = \frac{3}{m}\sum_{k=1}^m x_{ik}x_{jk}\, .
\end{equation}
The above formula shows that the sizes of the partial derivatives of $f$ remain stable as $n\ra\infty$. 
To establish the low complexity gradient condition, we need to show that there is a set of size $e^{o(n)}$ such that for any $x\in[0,1]^n$, $\nabla f(x)$ is within Euclidean distance $o(\sqrt{n})$ from one of the elements of this set. This is proved as a consequence of two key observations: 
\begin{enumerate}[1.]
\item By equation \eqref{tder},
\begin{align*}
\qquad \|\nabla f(x)-\nabla f(y)\|^2 &= \frac{9}{m^2}\sum_{i,j,k,l} (x_{ik}x_{jk}-y_{ik}y_{jk})(x_{il}x_{jl}-y_{il}y_{jl})\, .
\end{align*}
Expand the brackets on the right and consider one pair of terms in the expansion, say,
\begin{align*}
\frac{1}{m^2}\sum_{i,j,k,l} (x_{ik}x_{jk}x_{il}x_{jl}-x_{ik}x_{jk}y_{il}y_{jl})\, .
\end{align*}
This term may be written in a telescoping manner as 
\begin{align*}
\frac{1}{m^2}\sum_{i,j,k,l} x_{ik}x_{jk}x_{il}(x_{jl}-y_{jl}) + \frac{1}{m^2}\sum_{i,j,k,l} x_{ik}x_{jk}(x_{il}-y_{il})y_{jl}\, .
\end{align*}
Let $M(x)$ be the matrix whose $(i,j)^{\textup{th}}$ entry is $x_{ij}$. Consider the first sum. If $i$ and $k$ are fixed, then the sum in $j$ and $l$ is a quadratic form of the matrix $M(x)-M(y)$. This shows that the first sum is bounded by $m\|M(x)-M(y)\|_{\text{op}}$, 
where $\|A\|_{\text{op}}$ denotes the $L^2$ operator norm of a matrix $A$. Similarly bounding other terms, we get  
\begin{equation}\label{mainmatrix}
 \|\nabla f(x)-\nabla f(y)\|^2\le Cm\|M(x)-M(y)\|_{\text{op}}\, ,
\end{equation}
where $C$ is a universal constant. 
\item Let ${\mathcal M}$ denote the set of $m\times m$ symmetric matrices with entries in $[0,1]$. This set has low  complexity in operator norm, in the sense that it has a subset of size $e^{o(m^2)}$ such that for any $M\in {\mathcal M}$, there is some $N$ in this subset for which $\|M-N\|_{\textup{op}} = o(m)$. To see this, observe that if $\lambda_1,\ldots, \lambda_m$ are the eigenvalues of $M$, then $\sum_i \lambda_i^2 = \tr(M^2) \le m^2$. 
This implies that the $k^{\textup{th}}$-largest eigenvalue is bounded by $mk^{-1/2}$. Therefore there is a rank $k$ matrix $N$ such that $\|M-N\|_{\textup{op}} \le mk^{-1/2}$. Since a matrix of rank $k$ is determined by $k$ eigenvectors and eigenvalues, it is easy to see that the set of rank $k$ matrices has low complexity. Letting $k$ grow slowly with $m$, this allows us to establish the low complexity of ${\mathcal M}$. The proof of the low complexity of $\nabla f$ is completed by combining this information with the inequality \eqref{mainmatrix}. 
\end{enumerate}
\end{example}

\begin{example}[Three-term arithmetic progressions]\label{aex}
Index the elements of $[0,1]^n$ as $x = (x_i)_{i\in \zz/n\zz}$, and define $f:[0,1]^n\to \rr$  as 
\begin{align*}
f(x):= \frac{1}{n}\sum_{i,j\in \zz/n\zz} x_i x_{i+j} x_{i+2j}\, ,
\end{align*}
where the sums $i+j$ and $i+2j$ are carried out modulo $n$. Note that if $A$ is a subset of $\zz/n\zz$ and $x_i = 1$ if $i\in A$ and $0$ otherwise, then $f(x)$ is equal to the number of three-term arithmetic progressions in $A$ divided by $n$.  Now,
\[
\fpar{f}{x_i} = \frac{1}{n} \sum_{j\in \zz/n\zz} (x_{i+j} x_{i+2j} + x_{i-j} x_{i+j} + x_{i-2j} x_{i-j})\, .
\]
This shows that the sizes of the partial derivatives of $f$ remain stable as $n\ra\infty$. It turns out that this $f$, too, satisfies the low complexity gradient condition. The proof is not short enough to be presented here in full details (see \cite{chatterjeedembo14} for that), but the main idea may be easily explained as follows. 
The discrete Fourier transform $\hx$ of a point $x\in\rr^{\zz/n\zz}$ is defined as 
\[
\hx_j := \frac{1}{\sqrt{n}} \sum_k x_k e^{2\pi\mathrm{i} jk/n}\, ,
\]
where $\mathrm{i} = \sqrt{-1}$. The following inequality was proved in \cite{chatterjeedembo14}:
\[
\|\nabla f(x)-\nabla f(y)\|^2 \le Cn^{1/2}\max_i |\hx_i - \hat{y}_i|\, .
\]
It is possible to show using this inequality that if we know a few of the large Fourier coefficients of $x$, then the vector $\nabla f(x) = (f_1(x),\ldots,f_n(x))$  can be approximately recovered. This allows us to establish the low complexity gradient condition for $f$.
\end{example}

\section{Nonlinear large deviations}\label{nonlinearsec}
Having discussed several examples of functions satisfying the low complexity gradient condition, let us now review the main result of \cite{chatterjeedembo14},  which says that the approximation \eqref{star} holds for such functions. Let us begin with a limiting statement that applies when $p$ is fixed and $n$ tends to infinity. This is not pertinent to large deviations for sparse random graphs, but gives a nice, clean result that encapsulates the essence of the low complexity gradient condition. This result is powerful enough to imply the large deviation results for dense random graphs.

For each $n$, let $f_n:[0,1]^n \ra \rr$ be a twice differentiable function. For each $n$, $i$ and $j$, let $c_{n,i,j}$ be a uniform upper bound on the mixed partial derivative 
\[
\biggl|\mpar{f_n}{x_i}{x_j}\biggr|\,.
\]
Additionally, assume that $a$ and $b$ are constants such that for each $n$, $|f_n|$ is bounded by $an$ and the absolute values of all the first-order partial derivatives of $f_n$ are bounded by $b$. Let $\phi_{n,p}$ be defined according to the formula~\eqref{phipdef} applied to $f_n$, and assume that $\phi_p(t) := \lim_{n\ra \infty} \phi_{p,n}(t)$ 
exists for each $t$.
\begin{theorem}\label{easythm}
Consider the setting introduced above. Let $C(f_n)$ be the complexity of the gradient of $f_n$, as defined in Section \ref{lowsec}. Suppose that as $n\ra \infty$, the following conditions hold:
\begin{align}\label{easycond}
C(f_n) = o(1), \ \ \  \sum_{i,j} c_{n,i,j} = o(n^2), \ \ \ \sum_{i,j} c_{n,i,j}^2 = o(n),\ \ \ \sum_i c_{n,i,i} = o(n)\,.
\end{align}
Take any $p\in (0,1)$ and let $Y = (Y_1,\ldots, Y_n)$ be a vector of i.i.d.\ random variables with $\pp(Y_i)=p=1-\pp(Y_i=0)$. Let $\phi_p$ be defined as above. Then for any $t$ where $\phi_p$ is continuous and finite,
\[ 
\lim_{n\ra \infty}\frac{\log\pp(f_n(Y)\ge tn)}{n} = - \phi_p(t)\,.
\]
\end{theorem}
Revisiting the examples of Section \ref{lowsec}, recall that in three out of the four cases we verified or at least sketched that $C(f_n)\ra 0$ as $n\ra \infty$. It is easy to check in each of theses cases that the remaining three conditions in \eqref{easycond} are valid. For instance,~if 
\[
f_n(x) = \frac{1}{n}\sum_{1\le i<j\le n} x_i x_j\,,
\]
then $c_{n,i,i}=0$ and $c_{n,i,j} = 1/n$ for $i\ne j$. Therefore, \eqref{easycond} holds and so we can apply Theorem \ref{easythm} for this function.

Theorem \ref{easythm} is a straightforward corollary of the main theorem of \cite{chatterjeedembo14}, which is stated below. This theorem gives a quantitative error bound instead of a limiting result, which can be used to analyze situations where $p\ra 0$ as $n\ra \infty$. %The exact statement of the theorem goes as follows.

Fix some $n\ge 1$. Let $\|f\|$ denote the supremum norm of $f:[0,1]^n \ra\rr$. Suppose that
$f:[0,1]^n \ra\rr$ is twice continuously differentiable in $(0,1)^n$, such 
that $f$ and all its first and second order derivatives extend continuously to the boundary. 
For each $i$ and $j$, let
\[
f_i := \fpar{f}{x_i}\ \  \text{ and } \ \ 
f_{ij} := \mpar{f}{x_i}{x_j}. 
\]
Define 
\[
a := \|f\|, \ \ b_i := \|f_i\| \ \ \text{ and } \ \ c_{ij} := \|f_{ij}\|\, . 
\]
Given $\ep>0$, let $\dd(\ep)$ be a finite subset of $\rr^n$ such that for all $x\in \{0,1\}^n$, there exists $d = (d_1,\ldots, d_n)\in \dd(\ep)$ such that 
\begin{equation*}%\label{entropy}
\sum_{i=1}^n (f_i(x)- d_i)^2\le n \ep^2.
\end{equation*}
The following result gives an error bound for the approximation \eqref{star} in terms of the quantities $a$, $b_i$, $c_{ij}$ and the sizes of the sets $\dd(\ep)$. 
\begin{theorem}[\cite{chatterjeedembo14}]\label{nldp}
Take $f$ as  above, $p\in (0,1)$, $Y$  as in Theorem \ref{easythm} and $\phi_p$ as in \eqref{phidef}.
Then, for any $\delta>0$, $\ep>0$ and $t\in \rr$,
\begin{align*}
\log \pp(f(Y)\ge tn) &\le - \phi_p(t-\delta) n + \textup{complexity term} \\
&\qquad + \textup{smoothness term}\, ,
\end{align*}
where with $a$, $b$, $c_{ij}$, $\dd(\ep)$ defined above, 
\begin{align*}
\textup{complexity term} &:= 
\frac{1}{4}\Big(n\sum_{i=1}^n\beta_i^2\Big)^{1/2}\ep + 3n\ep +\log \biggl(\frac{4\phi_p(t)(\frac{1}{n}\sum_{i=1}^n b_i^2)^{1/2}}{\delta\ep}\biggr) \\
&\qquad + \log \biggl|\dd\biggl(\frac{\delta\ep}{4\phi_p(t)}\biggr)\biggr|\, , \ \text{ and}\\
%\end{align*}
%and
%\begin{align*}
\textup{smoothness term} &:= 4\biggl(\sum_{i=1}^n(\alpha \gamma_{ii} + \beta_i^2)+ \frac{1}{4}\sum_{i,j=1}^n \bigl(\alpha \gamma_{ij}^2+\beta_i \beta_j \gamma_{ij} + 4\beta_i \gamma_{ij}\bigr)\biggr)^{1/2} \\
&\qquad + \frac{1}{4}\Big(\sum_{i=1}^n \beta_i^2\Big)^{1/2}\Big(\sum_{i=1}^n \gamma_{ii}^2\Big)^{1/2} + 3\sum_{i=1}^n \gamma_{ii}+ \log 2\, ,
\end{align*}
where
\begin{align*}
\alpha &:= n\phi_p(t) + n|\log p| + n|\log (1-p)|\, ,\\
\beta_i &:= \frac{2\phi_p(t) b_i}{\delta} + |\log p| + |\log(1-p)|\, , \text{ and}\\
\gamma_{ij}  &:= \frac{2\phi_p(t) c_{ij}}{\delta} +\frac{ 6\phi_p(t) b_i b_j}{n\delta^2}\, .
\end{align*}
Moreover,  
\begin{align*}
\log \pp(f(Y)\ge tn) &\ge -\phi_p(t+\delta_0)n - \ep_0 n - \log 2\, ,
\end{align*}
where
\[
\ep_0  := \frac{1}{\sqrt{n}}\biggl(4+\biggl|\log\frac{p}{1-p}\biggr|\biggr)
\]
and
\[
\delta_0 := \frac{2}{n}\biggl(\sum_{i=1}^n (ac_{ii} + b_i^2)\biggr)^{1/2}\, .
\]
\end{theorem}

When applying Theorem \ref{nldp} to a given problem, one needs to first compute an error bound depending on some specific choice of $\ep$ and $\delta$, and then optimize the resulting bound over all possible values of $\ep$ and $\delta$. 
\begin{example}[Triangles in $G(n,p)$]
To get a flavor of the consequences of Theorem~\ref{nldp}, let us look at what it says for triangles in $G(n,p)$. Recall the function $I_p:[0,1]\ra \rr$ defined in~\eqref{ipdef}. For $x = (x_{ij})_{1\le i<j\le n}\in [0,1]^{n(n-1)/2}$, define
\[
I_p(x) := \sum_{1\le i<j\le n} I_p(x_{ij})
\]
and 
\[
T(x) := \frac{1}{6}\sum_{i,j,k=1}^n x_{ij}x_{jk}x_{ki}\, ,
\]
where $x_{ji}=x_{ij}$ and $x_{ii}=0$. For $u > 1$ define 
\[
\psi_p(u) := \inf\{I_p(x): T(x)\ge u\, \ee(T_{n,p})\}\, ,
\]
where $T_{n,p}$ is the number of triangles in $G(n,p)$. The following result was proved in~\cite{chatterjeedembo14}. The complexity calculations of Example \ref{tex} are crucial for this result. 
\begin{theorem}[\cite{chatterjeedembo14}]\label{cdthm}
For $u> 1$, $n$ sufficiently large (depending only on $u$), and $n^{-1/6}\le p\le 1-n^{-1}$, 
\begin{align*}
1-\frac{c \log n}{n^{1/6} p^{2}} \le \frac{\psi_p(u)}{-\log \pp(T_{n,p}\ge u\, \ee(T_{n,p}))} \le 1+ \frac{C (\log n)^{33/29}}{n^{1/29} p^{42/29}} \, ,
\end{align*}
where $c$ and $C$ are constants that depend only on $u$.
\end{theorem}
Theorem \ref{cdthm} is a non-asymptotic result. To get an asymptotic statement,  note that
\[
\frac{\psi_p(u)}{-\log \pp(T_{n,p}\ge u\, \ee(T_{n,p}))} \ra 1
\]
if $n\ra\infty$ and  $p\ra 0$ slower than  $n^{-1/42}(\log n)^{11/14}$. Theorem \ref{lzthm} was proved in~\cite{lubetzky12} by analyzing the asymptotic behavior of $\psi_p(u)$. 
\end{example}

\begin{example}[Three-term arithmetic progressions]\label{threeterm}
Fix $n \in \N$ and $p \in (0,1)$. Let $A$ be a random subset of $\zz/n\zz$, constructed by keeping each element with probability $p$, and dropping with probability $1-p$, independently of each other. Recall the function $I_p:[0,1]\ra \rr$ defined in equation \eqref{ipdef}. For $x = (x_i)_{i\in \zz/n\zz}$, let 
\[
I_p(x) := \sum_{i\in \zz/n\zz} I_p(x_i)\,.
\]
The following large deviation result about the number of three-term arithmetic progressions in $A$ was proved in \cite{chatterjeedembo14}. Again, the complexity calculations of Example~\ref{aex} are used for proving this theorem.
\begin{theorem}[\cite{chatterjeedembo14}]\label{arith}
Let $A$ be a random subset of $\zz/n\zz$, constructed as above. Let $X$ be the number of pairs $(i,j)\in (\zz/n\zz)^2$ such that $\{i,i+j, i+2j\}\subseteq A$. Let $I_p$ be defined as above. Let 
\begin{align*}
\theta_p(u) := \inf\biggl\{&I_p(x): x\in [0,1]^{\zz/n\zz} \ \textup{such that $\sum_{i,j\in \zz/n\zz} x_i x_{i+j} x_{i+2j}\ge u\, \ee(X)$}\biggr\}\, .
\end{align*}
Suppose that $n^{-1/162}\le p\le 1-n^{-1}$. Then for any $u>1$,
\begin{align*}
1 - c\, n^{-1/6}p^{-6}\log n&\le \frac{\theta_p(u)}{-\log \pp(X\ge u\, \ee(X))}  \le 1 + Cn^{-1/29}p^{-162/29}(\log n)^{33/29}\, ,
\end{align*}
where $C$ and $c$ are constants that depend only on $u$.
\end{theorem}
This theorem gives an approximation for the upper tail probabilities of the number of three-term arithmetic progressions in random subsets of $\zz/n\zz$ when $p$ is either fixed or decays slower than $n^{-1/162}(\log n)^{33/162}$ as $n\ra\infty$. Note that when $p=1/2$, calculating these upper tail probabilities is the same problem as counting the number of subsets of $\zz/n\zz$ that contain more than a given number of three-term progressions.

The study of arithmetic progressions in subsets of integers has a long and storied history, most of which is concerned with questions of existence. An excellent survey of old and new results is available in \cite{taovu}. Counting the number of sets with a given number of arithmetic progressions, or understanding the typical structure of sets that contain lots of progressions, are challenges of a different type, falling within the purview of large deviations theory. Recently a certain amount of interest has begun to grow around the resolution of such questions, quickly leading to the realization that conventional large deviations theory will not provide the answers. The most pertinent papers are the recent articles on probabilistic properties of the so-called `non-conventional averages' \cite{carinci, kifer, kv}. For example, \cite{carinci} gives a large deviation principle for sums of the type $\sum x_i x_{2i}$. 
\end{example}

\section{The naive mean field approximation}\label{meanfieldsec}
Theorem \ref{nldp} is  a consequence of a more general theorem about normalizing constants. In a nutshell, the theorem says that the so-called `naive mean field approximation' of statistical physics is valid when the low complexity gradient condition holds. The purpose of this section is to give the precise statement of this theorem and present a sketch of its proof.

Take any twice differentiable $f:[0,1]^n$ and let $a$, $b_i$, $c_{ij}$ and $\dd(\ep)$ be as in Theorem \ref{nldp}.
For $x = (x_1,\ldots, x_n)\in [0,1]^n$, let
\[
I(x) := \sum_{i=1}^n (x_i\log x_i + (1-x_i)\log (1-x_i))\, .
\]
Let 
\[
F:= \log \sum_{x\in [0,1]^n} e^{f(x)}\, .
\]
One version of the  naive mean field approximation is that 
\[
F\approx \sup_{x\in[0,1]^n}(f(x)-I(x))\,.
\]
The following theorem, proved in \cite{chatterjeedembo14}, gives a sufficient condition for the validity of this approximation. As far as I know, there is no other general sufficient condition for the validity of the naive mean field approximation.
\begin{theorem}[\cite{chatterjeedembo14}]\label{nldp2}
Let all notation be as above. Then for any $\ep>0$,
\begin{align*}
F &\le \sup_{x\in[0,1]^n}(f(x)-I(x)) + \frac{1}{4}\Big(n\sum_{i=1}^nb_i^2\Big)^{1/2}\ep + 3n\ep +\log |\dd(\ep)|\\
&\qquad + 4\biggl(\sum_{i=1}^n(ac_{ii} + b_i^2)+ \frac{1}{4}\sum_{i,j=1}^n \bigl(ac_{ij}^2+b_i b_j c_{ij} + 4b_i c_{ij}\bigr)\biggr)^{1/2} \\
&\qquad + \frac{1}{4}\Big(\sum_{i=1}^n b_i^2\Big)^{1/2}\Big(\sum_{i=1}^n c_{ii}^2\Big)^{1/2} + 3\sum_{i=1}^n c_{ii}+ \log 2\,,
\end{align*}
and 
\[
F \ge \sup_{x\in [0,1]^n} (f(x)-I(x)) - \frac{1}{2}\sum_{i=1}^n c_{ii}\, .
\]
\end{theorem}
As in Theorem \ref{nldp}, the error bound needs to optimized over $\ep$ when applying Theorem \ref{nldp2}. Theorem~\ref{nldp} is deduced from Theorem~\ref{nldp2} by replacing $f$ in Theorem~\ref{nldp} with a function $g$ which is a smooth approximation of the function that equals $1$ where $f(x)\ge tn$ and $0$ where $f(x)<tn$.

To finish the discussion, let us now see a sketch of the proof of Theorem \ref{nldp2}. This sketch is given with the intention of being helpful to someone who may try to solve the open problems listed in the next section. 
Let $X=(X_1,\ldots, X_n)$ be a random vector that has probability density proportional to $e^{f(x)}$ on $\{0,1\}^n$ with respect to the counting measure.  For each $i$, define a function $\hx_i:[0,1]^n \ra [0,1]$ as
\[
\hx_i(x) = \ee(X_i\mid X_j = x_j, \, 1\le j\le n, \ j\ne i)\,. 
\]
Let $\hx:[0,1]^n \ra[0,1]^n$ be the vector-valued function whose $i^{\textup{th}}$ coordinate function is $\hx_i$.  Let $\HX = \hx(X)$.  The first step in the proof is to show that if the smoothness term is small, then $f(X) \approx f(\HX)$ with high probability.
To show this, define 
\[
h(x) := f(x) - f(\hx(x))\, .
\]
Let $u_i(t, x) := f_i(tx +(1-t)\hx(x))$, 
so that
\[
h(x) = \int_0^1 \sum_{i=1}^n (x_i-\hx_i(x)) u_i(t,x)\, dt\,. 
\]
Thus, if $D := f(X)-f(\HX)$, then 
\begin{equation}\label{dagger}
\ee(D^2) = \int_0^1 \sum_{i=1}^n \ee((X_i-\HX_i) u_i(t, X) D) \, dt\, . 
\end{equation}
Let $X^{(i)}$ denote the random vector $(X_1,\ldots, X_{i-1}, 0, X_{i+1}, \ldots, X_n)$ and let $D_i := h(X^{(i)})$.  Then note that $ u_i(t, X^{(i)}) D_i$ is a function of the random variables $(X_j)_{j\ne i}$ only.  Therefore since $\HX_i = \ee(X_i\mid (X_j)_{j\ne i})$, 
\[
\ee((X_i-\HX_i) u_i(t, X^{(i)}) D_i) = 0\,. 
\]
Thus, 
\begin{align*}
&\ee((X_i-\HX_i) u_i(t, X) D) \\
&= \ee((X_i-\HX_i) u_i(t, X) D)  - \ee((X_i-\HX_i) u_i(t, X^{(i)}) D_i)\, . 
\end{align*}
If the smoothness term is small, then  $u_i(t, X) \approx u_i(t, X^{(i)})$ and $D \approx D_i$. Together with the identity \eqref{dagger}, 
this shows that $f(X)\approx f(\HX)$ with high probability. 

Next, define a function $g:[0,1]^n \times [0,1]^n \ra\rr$ as 
\[
g(x,y) := \sum_{i=1}^n (x_i \log y_i + (1-x_i)\log (1-y_i))\, .
\]
By a similar argument as above, it is possible to show that if the smoothness term is small, then with high probability,
\begin{equation*}\label{main2}
g(X, \HX) \approx g(\HX, \HX) = I(\HX)\, .
\end{equation*}
Let $A$ be the set of all $x$ where $f(x) \approx f(\hx(x))$ and $g(x, \hx(x)) \approx I(\hx(x))$.  Since $X\in A$ with high probability,
\[
\frac{\sum_{x\in A} e^{f(x)}}{\sum_{x\in \{0,1\}^n} e^{f(x)}} \approx 1\, .
\]
Therefore 
\begin{align*}
F &= \log \sum_{x\in \{0,1\}^n} e^{f(x)}\approx \log \sum_{x\in A} e^{f(x)}  \\
&\approx \log \sum_{x\in A} e^{f(\hx(x)) -I(\hx(x))+g(x, \hx(x))}\, .
\end{align*}
Now let $\ep$ be a small positive number.  Using the set $\dd(\ep)$, it is easy to produce a set $\dd'(\ep)\subseteq [0,1]^n$ such that $|\dd(\ep)|=|\dd'(\ep)|$, and  for each $x$ there exists $p\in \dd'(\ep)$ such that $\hx(x)\approx p$.  
 For each $p\in \dd'(\ep)$ let $\mathcal{P}(p)$ be the set of all $x\in \{0,1\}^n$ such that $\hx(x)\approx p$. The crucial fact is that for any $p\in [0,1]^n$, 
\[
\sum_{x\in \{0,1\}^n} e^{g(x,p)} = 1\,. 
\]
Therefore,
\begin{align*}
&\log \sum_{x\in A} e^{f(\hx(x))-I(\hx(x)) + g(x, \hx(x) )} \\
&\le \log \sum_{p\in \dd'(\ep)}\sum_{x\in \mathcal{P}(p)} e^{f(\hx(x))-I(\hx(x)) + g(x, \hx(x))}  \\
&\approx \log \sum_{p\in \dd'(\ep)}\sum_{x\in \mathcal{P}(p)} e^{f(p)-I(p) + g(x,p)}\\
&\le    \log \sum_{p\in \dd'(\ep)} e^{f(p)-I(p)}\label{fineq3o}
 \le \log|\dd'(\ep)| + \sup_{p\in [0,1]^n} (f(p)-I(p))\, .
\end{align*}
This completes the proof sketch for the upper bound. The lower bound is much more straightforward. Take any $y\in [0,1]^n$. Let $Y = (Y_1,\ldots, Y_n)$ be independent $0$-$1$ random variables, with $\pp(Y_i=1)=y_i$. Then by Jensen's inequality,
\begin{align*}
\sum_{x\in \{0,1\}^n} e^{f(x)} &= \sum_{x\in \{0,1\}^n} e^{f(x)-g(x,y)+g(x,y)}\\
&= \ee(e^{f(Y)-g(Y,y)})\\
&\ge \exp(\ee(f(Y) - g(Y,y)))\\
&= \exp(\ee(f(Y)) - I(y))\, .
\end{align*}
Then, by using arguments similar to those employed for the upper bound, one can prove that if the error term in the lower bound is small, then $\ee(f(Y)) \approx f(y)$. Since this is true for any $y$, this completes the sketch of the proof of the lower bound.

\section{Open problems}
There are many open problems about large deviations for random graphs, since the subject is still in a stage of development. The following is a partial list of the most important questions.
\begin{enumerate}[1.]
\item Produce explicit non-constant solutions of the variational problems arising from applications of Theorem \ref{cvthm}. Currently, only the existence of non-constant solutions is known in certain regimes, but explicit non-constant solutions, or mathematically provable qualitative properties of non-constant solutions, are unknown. This is important for understanding the conditional behavior of dense random graphs if some rare event takes place.
\item As a concrete instance of the above general question, analyze the behavior of the non-constant solutions of \eqref{phidef}.
\item Improve the main nonlinear large deviation result (Theorem \ref{nldp}), so that results like Theorem \ref{lzthm} can be proved when $p$ tends to zero at an optimal rate.
\item As an example of the above, show that Theorem \ref{lzthm} holds when $p\ra 0$ slower than $n^{-1/2}$.
\item Develop a sparse regularity lemma and a sparse graph limit theory that is powerful enough to prove results like Theorem \ref{lzthm} and Theorem \ref{bthm}. In fact, a reasonable test for the completeness of a sparse graph limit theory is whether it can lead to a solution of the large deviation question for sparse Erd\H{o}s--R\'enyi random graphs. This is because analyzing the large deviation behavior of $G(n,p)$ for small $p$ requires a full understanding of {\it all} possible sparse graph structures rather than focusing a small subset of graphs with nice properties.
\item Extend the large deviation results for three-term arithmetic progressions (Example \ref{threeterm}) to longer progressions. The reader may recall that discrete Fourier analysis was used in the analysis of three-term progressions. The method does not seem to extend easily to longer progressions. It is possible that higher order Fourier analysis (Gowers norms) or the hypergraph regularity lemma may be needed for longer progressions.
\item Find explicit solutions to the variational problems coming from arithmetic progressions, in the spirit of Theorems \ref{triangle}, \ref{lzthm} and \ref{bthm} of this survey.
\item In addition to the above problems, there are many related open problems of similar flavor about exponential random graphs. In particular, non-constant solutions of the variational problem \eqref{expvar} are generally not known, except for the solution of a related problem in \cite{kenyon14}. Sparse exponential random graphs are still out of the reach of mathematical results, except for some progress in \cite{chatterjeedembo14}.
\end{enumerate}

\section*{Acknowledgments}
I thank Amir Dembo, Jafar Jafarov, Jun Yan and Yufei Zhao for carefully reading the first draft and giving useful comments, and Rick Durrett for suggesting that I write this survey.

\bibliographystyle{amsplain}

\end{document}